\documentclass[12pt]{amsart}
\usepackage{hyperref}

\newtheorem{theorem}{Theorem}[section]
\newtheorem{definition}[theorem]{Definition}

\newtheorem{proposition}[theorem]{Proposition}

\newtheorem{question}[theorem]{Question}

\begin{document}

\title{Quantum isometries and noncommutative spheres}

\author{Teodor Banica}
\address{T.B.: Department of Mathematics, Toulouse 3 University, 118 route de Narbonne, 31062 Toulouse, France. {\tt banica@math.ups-tlse.fr}}

\author{Debashish Goswami}
\address{D.G.: Theoretical Statistics and Mathematics Unit, 203 Barrackpore Trunk Road, Kolkata 700 108, India. {\tt goswamid@isical.ac.in}}


\begin{abstract}
We introduce and study two new examples of noncommutative spheres: the half-liberated sphere, and the free sphere. Together with the usual sphere, these two spheres have the property that the corresponding quantum isometry group is ``easy'', in the representation theory sense. We present as well some general comments on the axiomatization problem, and on the ``untwisted'' and ``non-easy'' case.
\end{abstract}

\maketitle

\section*{Introduction}

The aim of the present paper is to bring some contributions to the theory of noncommutative spheres, by using a number of ideas and tools coming from the recent work on quantum isometry groups \cite{bg2}, \cite{gos}, and on easy quantum groups \cite{bsp}, \cite{bve}.

The noncommutative spheres were introduced by Podle\'s in \cite{pod}, as twists of the usual spheres. The natural framework for the study of such noncommutative objects is Connes' noncommutative geometry \cite{con}. This has led to a systematic study of the associated spectral triples, with the explicit computation of a number of related Riemannian geometric invariants. See Connes and Dubois-Violette \cite{cd1}, \cite{cd2}, Connes and Landi \cite{cla}, Dabrowski, D'Andrea, Landi and Wagner \cite{dd+}.

A useful, alternative point of view comes from the relationship with the quantum groups. The structure of the usual sphere $S^{n-1}$ is intimately related to that of the orthogonal group $O_n$, and when twisting the sphere the orthogonal group gets twisted as well, and becomes a quantum group. See Varilly \cite{var}.

The recently developed theory of quantum isometry groups \cite{bg2}, \cite{gos} provides a good abstract framework for the study of the exact relationship between noncommutative spheres and quantum groups. Some key preliminary results in this direction were obtained in \cite{bg1}, where it was shown that the usual spheres do not have quantum symmetry, and in \cite{bg3}, where the case of the Podle\'s sphere is studied in detail.

In this paper we present some more general results in this direction, mixing spectral triple and quantum group techniques. The idea is that the ``easy'' (and ``untwisted'') orthogonal quantum groups, introduced in \cite{bsp} and classified in \cite{bve}, are as follows: the orthogonal group $O_n$, the half-liberated orthogonal group $O_n^*$, and the free orthogonal group $O_n^+$. Together with the above considerations, this suggests that in the ``untwisted'' and ``easy'' case we should have exactly 3 examples of noncommutative spheres: the usual sphere $S^{n-1}$, a half-liberated sphere $S^{n-1}_*$, and a free sphere $S^{n-1}_+$. We will present here a number of results in this direction, notably by introducing and studying in detail the half-liberated sphere $S^{n-1}_*$, and the free sphere $S^{n-1}_+$.

The paper is organized as follows: in 1-2 we discuss the construction of the three spheres, in 3-5 we study the associated projective spaces and spherical integrals, and in 6-7 we work out the relation with the quantum isometry groups. 

The final section, 8, contains a few concluding remarks.

\subsection*{Acknowledgements}

We would like to thank J. Bichon and S. Curran for several useful discussions. T.B. was supported by the ANR grants ``Galoisint'' and ``Granma'', and D.G. was supported by the project ``Noncommutative Geometry and Quantum Groups'', funded by the Indian National Science Academy.

\section{The usual sphere}

The simplest example of noncommutative sphere is the usual sphere $S^{n-1}\subset\mathbb R^n$. Our first goal will be to find a convenient functional analytic description of it. In this section we present a series of 5 basic statements in this direction, all well-known, and all to be extended in the next section to the noncommutative setting.

\begin{theorem}
The sphere $S^{n-1}\subset\mathbb R^n$ is the spectrum of the $C^*$-algebra
$$A_n=C^*\left(x_1,\ldots,x_n\Big|x_i=x_i^*, x_ix_j=x_jx_i,\sum x_i^2=1\right)$$
generated by $n$ self-adjoint commuting variables, whose squares sum up to $1$.
\end{theorem}

\begin{proof}
The first remark is that the algebra in the statement is indeed well-defined, due to the condition $\Sigma\,x_i^2=1$, which shows that we have $||x_i||\leq 1$ for any $C^*$-norm.

Consider the algebra $A_n'=C(S^{n-1})$, with the standard coordinates denoted $x_i'$. By the universal property of $A_n$ we have a morphism $A_n\to A_n'$ mapping $x_i\to x_i'$.

On the other hand, by the Gelfand theorem we have $A_n=C(X)$ for a certain compact space $X$, and we can define a map $X\to S^{n-1}$ by $p\to (x_1(p),\ldots,x_n(p))$. By transposing we get a morphism $A_n'\to A_n$ mapping $x_i'\to x_i$, and we are done.
\end{proof}

The second ingredient that we will need is a functional analytic description of the action of the orthogonal group $O_n$ on $S^{n-1}$. This action can be seen as follows.

\begin{theorem}
We have a coaction $\Phi:A_n\to C(O_n)\otimes A_n$ given by
$$\Phi(x_i)=\sum_ju_{ij}\otimes x_j$$
where $u_{ij}\in C(O_n)$ are the standard matrix coordinates: $u_{ij}(g)=g_{ij}$.
\end{theorem}

\begin{proof}
Consider indeed the action map $O_n\times S^{n-1}\to S^{n-1}$, given by $(g,p)\to gp$. By transposing we get a coaction map $\Phi$ as in the statement.
\end{proof}

The uniform measure on $S^{n-1}$ is the unique probability measure which is invariant under the action of $O_n$. In functional analytic terms, the result is as follows.

\begin{theorem}
There is a unique positive unital trace $tr:A_n\to\mathbb C$ satisfying the invariance condition $(id\otimes tr)\Phi(x)=tr(x)1$.
\end{theorem}

\begin{proof}
We can define indeed $tr:A_n\to\mathbb C$ to be the integration with respect to the uniform measure on $S^{n-1}$: the positivity condition follows from definitions, and the invariance condition as in the statement follows from $dp=d(gp)$, for any $g\in O_n$.

Conversely, it follows from the general theory of the Gelfand correspondence that a trace as in the statement must come from the integration with respect to a probability measure on $S^{n-1}$ which is invariant under $O_n$, and this gives the uniqueness.
\end{proof}

\begin{theorem}
The canonical trace $tr:A_n\to\mathbb C$ is faithful.
\end{theorem}

\begin{proof}
This follows from the well-known fact that the uniform measure on the sphere takes a nonzero value on any open set. 
\end{proof}

Finally, we have the well-known result stating that $S^{n-1}$ can be identified with the first slice of $O_n$. In functional analytic terms, the result is as follows.

\begin{theorem}
The following algebras, with generators and traces, are isomorphic:
\begin{enumerate}
\item The algebra $A_n$, with generators $x_1,\ldots,x_n$, and with the trace functional.

\item The algebra $B_n\subset C(O_n)$ generated by $u_{11},\ldots,u_{1n}$, with the integration.
\end{enumerate}
\end{theorem}

\begin{proof}
From the universal property of $A_n$ we get a morphism $\pi:A_n\to B_n$ mapping $x_i\to u_{1i}$. The invariance property of the integration functional $I:C(O_n)\to\mathbb C$ shows that $tr'=I\pi$ satisfies the invariance condition in Theorem 1.3, so we have $tr=tr'$. Finally, from the faithfulness of $tr$ we get that $\pi$ is an isomorphism, and we are done.
\end{proof}

\section{Noncommutative spheres}

We are now in position of introducing two basic examples of noncommutative spheres: the half-liberated sphere, and the free sphere. The idea will be of course to weaken or simply remove the commutativity conditions in Theorem 1.1.

\begin{definition}
We consider the universal $C^*$-algebras
\begin{eqnarray*}
A_n^*&=&C^*\left(x_1,\ldots,x_n\Big|x_i=x_i^*,x_ix_jx_k=x_kx_jx_i,\sum x_i^2=1\right)\\
A_n^+&=&C^*\left(x_1,\ldots,x_n\Big|x_i=x_i^*,\sum x_i^2=1\right)
\end{eqnarray*}
generated by $n$ self-adjoint variables whose squares sum up to $1$, subject to the half-commutation relations $abc=bca$, and to no relations at all.
\end{definition}

Our next goal will be to find suitable analogues of Theorems 1.2, 1.3, 1.4, 1.5. For this purpose, we will need appropriate ``noncommutative versions'' of the orthogonal group $O_n$. So, let us consider the following universal algebras:
\begin{eqnarray*}
C(O_n^*)&=&C^*\left(u_{11},\ldots, u_{nn}\Big|u_{ij}=u_{ij}^*,u_{ij}u_{kl}u_{st}=u_{st}u_{kl}u_{ij},u^t=u^{-1}\right)\\
C(O_n^+)&=&C^*\left(u_{11},\ldots,u_{nn}\Big|u_{ij}=u_{ij}^*,u^t=u^{-1}\right)
\end{eqnarray*}

These algebras, introduced in \cite{bsp}, \cite{wan}, are Hopf algebras in the sense of Woronowicz \cite{wo1}, \cite{wo2}. We refer to the recent paper \cite{bve} for a full discussion here.

With these definitions in hand, we can state and prove now an analogue of Theorem 1.2. We agree to use the generic notation $A_n^\times$ for the 3 algebras constructed so far.

\begin{theorem}
We have a coaction $\Phi:A_n^\times\to C(O_n^\times)\otimes A_n^\times$ given by
$$\Phi(x_i)=\sum_ju_{ij}\otimes x_j$$
where $u_{ij}\in C(O_n^\times)$ are the standard generators.
\end{theorem}

\begin{proof}
We have to construct three maps, and prove that they are coactions. We will deal with all 3 cases at the same time. Consider the following elements:
$$X_i=\sum_ju_{ij}\otimes x_j$$

These elements are self-adjoint, and their squares sum up to $1$. Moreover, in the case where both sets $\{x_i\}$ and $\{u_{ij}\}$ consist of commuting or half-commuting elements, the set $\{X_i\}$ consists by construction of commuting or half-commuting elements. 

These observations give the existence of maps $\Phi$ as in the statement. The fact that these maps satisfy the condition $(\Delta\otimes id)\Phi=(\Phi\otimes id)\Phi$ is clear from definitions.
\end{proof}

\begin{theorem}
There is a unique positive unital trace $tr:A_n^\times\to\mathbb C$ satisfying the condition $(id\otimes tr)\Phi(x)=tr(x)1$.
\end{theorem}

\begin{proof}
Consider the algebra $B_n^\times\subset C(O_n^\times)$ generated by the elements $u_{11},\ldots,u_{1n}$. By the universal property of $A_n^\times$ we have a morphism $\pi:A_n^\times\to B_n^\times$ mapping $x_i\to u_{1i}$, and by composing with the restriction of the Haar functional $I:C(O_n^\times)\to\mathbb C$, we obtain a trace satisfying the invariance condition in the statement.

As for the uniqueness part, this is a quite subtle statement, which will ultimately come from a certain ``easiness'' property of our 3 spheres. Our first claim is that we have the following key formula, where $tr$ is the trace that we have just constructed:
$$(I\otimes id)\Phi=tr(.)1$$

In order to prove this claim, we use the orthogonal Weingarten formula \cite{bcs}, \cite{bc1}, \cite{bsp}. To any integer $k$ let us associate the following sets:
\begin{enumerate}
\item $D_k$: all pairings of $\{1,\ldots,k\}$.

\item $D_k^*$: all pairings with an even number of crossings of $\{1,\ldots,k\}$.

\item $D_k^+$: all noncrossing pairings of $\{1,\ldots,k\}$.
\end{enumerate}

The Weingarten formula tells us that the Haar integration over $O_n^\times$ is given by the following sum, where $W_{kn}=G_{kn}^{-1}$, with $G_{kn}(p,q)=n^{loops(p\vee q)}$, and where the $\delta$ symbols are $1$ if all the strings join pairs of equal indices, and $0$ if not:
$$I(u_{i_1j_1}\ldots u_{i_kj_k})=\sum_{p,q\in D_k^\times}\delta_p(i)\delta_q(j)W_{kn}(p,q)$$

So, let us go back now to our claim. By linearity it is enough to check the equality on a product of basic generators $x_{i_1}\ldots x_{i_k}$. The left term is as follows:
\begin{eqnarray*}
(I\otimes id)\Phi(x_{i_1}\ldots x_{i_k})
&=&\sum_{j_1\ldots j_k}I(u_{i_1j_1}\ldots u_{i_kj_k})x_{j_1}\ldots x_{j_k}\\
&=&\sum_{j_1\ldots j_k}\sum_{p,q\in D_k^\times}\delta_p(i)\delta_q(j)W_{kn}(p,q)x_{j_1}\ldots x_{j_k}\\
&=&\sum_{p,q\in D_k^\times}\delta_p(i)W_{kn}(p,q)\sum_{j_1\ldots j_k}\delta_q(j)x_{j_1}\ldots x_{j_k}
\end{eqnarray*}

Let us look now at the last sum on the right. In the free case we have to sum quantities of type $x_{j_1}\ldots x_{j_k}$, over all choices of multi-indices $j$ which fit into our given noncrossing pairing $q$, and just by using the condition $\Sigma x_i^2=1$, we conclude that the sum is 1. The same happens in the classical case, with the changes that our pairing $q$ can now be crossing, but we can use now the commutation relations $x_ix_j=x_jx_i$. Finally, the same happens as well in the half-liberated case, because the fact that our pairing $q$ has now an even number of crossings allows us to use the half-commutation relations $x_ix_jx_k=x_kx_jx_i$, in order to conclude that the sum to be computed is 1.

Summarizing, in all cases the sum on the right is 1, so we get:
\begin{eqnarray*}
(I\otimes id)\Phi(x_{i_1}\ldots x_{i_k})
&=&\sum_{p,q\in D_k^\times}\delta_p(i)W_{kn}(p,q)1
\end{eqnarray*}

On the other hand, another application of the Weingarten formula gives:
\begin{eqnarray*}
tr(x_{i_1}\ldots x_{i_k})1
&=&I(u_{1i_1}\ldots u_{1i_k})1\\
&=&\sum_{p,q\in D_k^\times}\delta_p(1)\delta_q(i)W_{kn}(p,q)1\\
&=&\sum_{p,q\in D_k^\times}\delta_q(i)W_{kn}(p,q)1
\end{eqnarray*}

Since the Weingarten function is symmetric in $p,q$, this finishes the proof of our claim. So, let us get back now to the original question. Let $\tau:A_n^\times\to\mathbb C$ be a trace satisfying the invariance condition in the statement. We have:
\begin{eqnarray*}
\tau (I\otimes id)\Phi(x)
&=&(I\otimes\tau)\Phi(x)\\
&=&I(id\otimes\tau)\Phi(x)\\
&=&I(\tau(x)1)\\
&=&\tau(x)
\end{eqnarray*}

On the other hand, according to our above claim, we have as well:
\begin{eqnarray*}
\tau (I\otimes id)\Phi(x)
&=&\tau(tr(x)1)\\
&=&tr(x)
\end{eqnarray*}

Thus we get $\tau=tr$, which finishes the proof.
\end{proof}

We do not have an analogue of Theorem 1.4, and best is to proceed as follows.

\begin{definition}
We agree to replace from now on $A_n^\times$ with its GNS completion with respect to the canonical trace $tr:A_n^\times\to\mathbb C$. 
\end{definition}

We actually believe that the canonical trace is faithful on the algebraic part, so that our replacement is basically not needed, but we don't have a proof for this fact.

\begin{theorem}
The following algebras, with generators and traces, are isomorphic:
\begin{enumerate}
\item The algebra $A_n^\times$, with generators $x_1,\ldots,x_n$, and with the trace functional.

\item The algebra $B_n^\times\subset C(O_n^\times)$ generated by $u_{11},\ldots,u_{1n}$, with the integration.
\end{enumerate}
\end{theorem}

\begin{proof}
Consider the map $\pi:A_n^\times\to B_n^\times$, already constructed in the proof of Theorem 2.3. The invariance property of the integration functional $I:C(O_n^\times)\to\mathbb C$ shows that $tr'=I\pi$ satisfies the invariance condition in Theorem 2.3, so we have $tr=tr'$. Together with the positivity of $tr$ and with the basic properties of the GNS construction, this shows that $\pi$ is an isomorphism, and we are done.
\end{proof}

\section{Projective spaces}

In this section we study the projective spaces associated to our noncommutative spheres. Let us first recall that the projective space over a field $\mathbb F$ is by definition $\mathbb F^n-\{0\}/\sim$, where $x\sim y$ when $y=\lambda x$ for some $\lambda\in\mathbb F$. We will use the notation $P^{n-1}$ for the real projective space, and $P^{n-1}_c$ for the complex projective space.

Let us introduce the following definition.

\begin{definition}
We denote by $C_n^\times$ the subalgebra $<x_ix_j>\subset A_n^\times$, taken together with the restriction of the canonical trace.
\end{definition}

The noncommutative projective space that we are interested in is by definition the spectrum of $C_n^\times$, viewed as a noncommutative compact measured space.

As a first remark, in the classical case we get indeed the real projective space.

\begin{theorem}
We have $C_n=C(P^{n-1})$.
\end{theorem}

\begin{proof}
First, since each product of coordinates $x_ix_j:S^{n-1}\to\mathbb R$ takes equal values on $p$ and $-p$, this product can be regarded as being a function on $P^{n-1}$. Now since the collection of functions $\{x_ix_j\}$ separates the points of $P^{n-1}$, we get the result.
\end{proof}

Quite surprisingly, in the half-liberated case we get the complex projective space. 

\begin{theorem}
We have $C_n^*=C(P^{n-1}_c)$.
\end{theorem}

\begin{proof}
First, the half-commutation relations $abc=cba$ give $abcd=cbad=cdab$ for any $a,b,c,d\in\{x_1,\ldots,x_n\}$, so the elements $x_ix_j$ commute indeed with each other. We have to prove that the Gelfand spectrum of $C_n^*=<x_ix_j>$ is isomorphic to $P^{n-1}_c$.

For this purpose, we use the isomorphism $PO_n^*\simeq PU_n$ established in \cite{bve}. This isomorphism is given by $u_{ij}u_{kl}\to v_{ij}v_{kl}^*$, where $u,v$ denote respectively the fundamental corepresentations of $O_n^*,U_n$. By restricting attention to the first row of coordinates, this gives an embedding $C_n^*\subset C(PU_n)$, mapping $x_ix_j\to v_{1i}v_{1j}^*$.

Consider now the complex sphere $S_c^{n-1}\subset\mathbb C^n$, with coordinates denoted $z_1,\ldots,z_n$. By performing the standard identification $v_{1i}=z_i$, coming from the unitary version of Theorem 1.5, we obtain an embedding $C_n^*\subset C(S_c^{n-1})$, mapping $x_ix_j\to z_i\bar{z}_j$.

The image of this embedding is the subalgebra of $C(S^{n-1}_c)$ generated by the functions $z_i\bar{z}_j$. By using the same argument as in the real case, this gives the result.
\end{proof}

In view of the above results, it is tempting to conjecture that the ``threefold way'' that we are currently developing for quantum groups, noncommutative spheres and noncommutative projective spaces is actually part of the usual real/complex/quaternionic ``threefold way'', originally discovered by Frobenius, and known to play a fundamental role in mathematical physics, according to Dyson's paper \cite{dys}. 

We have here the following question.

\begin{question}
Do we have $C_n^+=C(P^{n-1}_k)$?
\end{question}

This is of course a quite vague question. The symbol $P^{n-1}_k$ on the right is supposed to correspond to some kind of tricky ``quaternionic projective space''. Note that the space $P^{n-1}_k=\mathbb K^n-\{0\}/\sim$ appearing in the existing literature won't be suitable for our purposes, simply because the algebra $C_n^+$ on the left is noncommutative. What we would need is rather a ``twist'', in the spirit of the quantum projective spaces in \cite{dla}.

The main problem here is to construct a representation of $C_n^+$, by using the Pauli matrices. With a bit of luck, this representation can be shown to be faithful, and the corresponding result can be interpreted as answering the above question.

A first piece of evidence comes from \cite{bc2}, where a certain faithful representation of $C(S_4^+)$ is constructed, by using the Pauli matrices. This is probably quite different from what we need, but the main technical fact, namely that ``the combinatorics of the noncrossing partitions can be implemented by the Pauli matrices'', is already there.

A second piece of evidence comes from the results in \cite{bve}, which suggest that the free quantum groups might be actually supergroups. Once again, this kind of argument is quite speculative, and maybe a bit far away from the present considerations.

Let us end this section by recording a few modest facts about $C_n^+$.

\begin{proposition}
The algebras $C_n^+$ are as follows:
\begin{enumerate}
\item At $n=2$ we have $C_2^+=C_2^*$.

\item At $n\geq 3$ we have $C_n^+\neq C_n^*$.
\end{enumerate}
\end{proposition} 

\begin{proof}
(1) This follows either from the isomorphism $O_2^+=O_2^*$ established in \cite{bve}, or directly from definitions, by using the fact that $C_2^+$ is commutative.

(2) It is enough here to prove that $C_3^+$ is not commutative. For this purpose, we will use the positive matrices in $M_2(\mathbb C)$. These are matrices of the following form:
$$Y=\begin{pmatrix}p&a\\ \bar{a}&q\end{pmatrix}$$

Here $p,q\in\mathbb R$ and $a\in\mathbb C$ must be chosen such that both eigenvalues are positive, and this happens for instance when $p,q>0$ and $a\in\mathbb C$ is small enough.

Let us fix some numbers $p_i,q_i>0$ for $i=1,2,3$, satisfying $\Sigma p_i=\Sigma q_i=1$. For any choice of small complex numbers $a_i\in\mathbb C$ satisfying $\Sigma a_i=0$, the corresponding elements $Y_i$ constructed as above will be positive, and will sum up to $1$. Moreover, by carefully choosing the $a_i$'s, we can arrange as for $Y_1,Y_2,Y_3$ not to pairwise commute.

Consider now the matrices $X_i=\sqrt{Y_i}$. These are all self-adjoint, and their squares sum up to 1, so we get a representation $A_3^+\to M_2(\mathbb C)$ mapping $x_i\to X_i$. Now this representation restricts to a representation $C_3^+\to M_2(\mathbb C)$ mapping $x_i^2\to Y_i$, and since the $Y_i$'s don't commute, it follows that $C_3^+$ is not commutative, and we are done. 
\end{proof}

The above result suggests the following extra question regarding $C_n^+$: what is the Gelfand spectrum of the algebra $C_n^+/I$, where $I\subset C_n^+$ is the commutator ideal? 

Observe that the canonical arrow $C_n^+\to C_n^*$ and Theorem 3.3 tell us that this Gelfand spectrum must contain $P^{n-1}_c$. Moreover, Proposition 3.5 shows that at $n=2$ this inclusion is an equality. However, at $n=3$ already the answer is not clear.

\section{Probabilistic aspects}

We know from the previous sections that we have three basic examples of ``noncommutative spheres'', namely those corresponding to the algebras $A_n,A_n^*,A_n^+$. In this section and in the next one we investigate the key problem of computing the integral over these noncommutative spheres of polynomial quantities of type $x_{i_1}\ldots x_{i_k}$.

\begin{definition}
The polynomial spherical integrals will be denoted 
$$I=\int_{S^{n-1}_\times}x_{i_1}\ldots x_{i_k}\,dx$$
with this quantity standing for the complex number obtained as image of the well-defined element $x_{i_1}\ldots x_{i_k}\in A_n^\times$ by the well-defined trace functional $tr:A_n^\times\to\mathbb C$.
\end{definition}

The problem of computing such integrals has been heavily investigated in the last years, and a number of results are available from \cite{bc1}, \cite{bcz}, \cite{csn}, \cite{dif}. In what follows we will make a brief presentation of this material, by focusing of course to the applications to $S^{n-1}_\times$. We will present as well some new results, in the half-liberated case.

Let us begin our study with an elementary result.

\begin{proposition}
We have the formula
$$\int_{S^{n-1}_\times}x_{i_1}\ldots x_{i_k}\,dx=0$$
unless each $x_i$ appears an even number of times.
\end{proposition}

\begin{proof}
This follows from the fact that for any $i$ we have an automorphism of $A_n^\times$ given by $x_i\to -x_i$. Indeed, this automorphism must preserve the trace, so if $x_i$ appears an odd number of times, the integral in the statement satisfies $I=-I$, so $I=0$.
\end{proof}

The basic tool for computing spherical integrals is the Weingarten formula. Let us recall from section 2 that associated to any integer $k$ are the following sets:
\begin{enumerate}
\item $D_k$: all pairings of $\{1,\ldots,k\}$.

\item $D_k^*$: all pairings with an even number of crossings of $\{1,\ldots,k\}$.

\item $D_k^+$: all noncrossing pairings of $\{1,\ldots,k\}$.
\end{enumerate}

These sets can be regarded as being associated to our spheres $S^{n-1}_\times$, because they come from the representation theory of the associated quantum groups $O_n^\times$.

\begin{theorem}
We have the Weingarten formula
$$\int_{S^{n-1}_\times}x_{i_1}\ldots x_{i_k}\,dx=\sum_{p,q\in D_k^\times}\delta_p(i)W_{kn}(p,q)$$
where $W_{kn}=G_{kn}^{-1}$, with $G_{kn}(p,q)=n^{loops(p\vee q)}$, and where the $\delta$ symbol is $1$ if all the strings of $p$ join pairs of equal indices of $i=(i_1,\ldots,i_k)$, and is $0$ if not.
\end{theorem}

\begin{proof}
This follows from the Weingarten formula in \cite{csn}, \cite{bc1}, \cite{bcs}, via the identification in Theorem 2.5, and from the fact that the Weingarten matrix is symmetric in $p,q$.
\end{proof}

As a first application, we have the following result.

\begin{theorem}
With $n\to\infty$, the standard coordinates of $S^{n-1}_\times$ are as follows:
\begin{enumerate}
\item Classical case: real Gaussian, independent.

\item Half-liberated case: symmetrized Rayleigh, their squares being independent.

\item Free case: semicircular, free.
\end{enumerate}
\end{theorem}

\begin{proof}
This follows from Theorem 4.3 and from the fact that $W_{kn}$ is asymptotically diagonal, see \cite{bc1}, \cite{bsp}. The only new assertion is the independence one in (2), which can be proved as in \cite{bsp}, by using the fact that the mixed cumulants vanish. 

Note that the independence in (1,2) follows as well from the exact formulae in Theorems 5.1 and 5.2 below, by letting $n\to\infty$ and by using the Stirling formula.
\end{proof}

\section{Spherical integrals}

We discuss in this section a quite subtle problem, of theoretical physics flavor, namely the exact computation of the polynomial integrals over $S^{n-1}_\times$. 

In the classical case, we have the following well-known result.

\begin{theorem}
The spherical integral of $x_{i_1}\ldots x_{i_k}$ vanishes, unless each $a\in\{1,\ldots,n\}$ appears an even number of times in the sequence $i_1,\ldots,i_k$. If $l_a$ denotes this number of occurrences, then
$$\int_{S^{n-1}}x_{i_1}\ldots x_{i_k}\,dx=\frac{(n-1)!!l_1!!\ldots l_n!!}{(n+\Sigma l_i-1)!!}$$
with the notation $m!!=(m-1)(m-1)(m-5)\ldots$
\end{theorem}

\begin{proof}
The first assertion follows from Proposition 4.2. The second assertion is well-known, and can be proved by using spherical coordinates, the Fubini theorem, and some standard partial integration tricks. See e.g. \cite{bcs}.
\end{proof}

In the case of the half-liberated sphere, we have the following result.

\begin{theorem}
The half-liberated spherical integral of $x_{i_1}\ldots x_{i_k}$ vanishes, unless each number $a\in\{1,\ldots,n\}$ appears the same number of times at odd and at even positions in the sequence $i_1,\ldots,i_k$. If $l_a$ denotes this number of occurrences, then:
$$\int_{S^{n-1}_*}x_{i_1}\ldots x_{i_k}\,dx=4^{\Sigma l_i}\frac{(2n-1)!l_1!\ldots l_n!}{(2n+\Sigma l_i-1)!}$$
\end{theorem}

\begin{proof}
First, by using Proposition 4.2 we see that the integral $I$ in the statement vanishes, unless $k=2l$ is even. So, assume that we are in the non-vanishing case. By using Theorem 3.3 the corresponding integral over the complex projective space $P^{n-1}_c$ can be viewed as an integral over the complex sphere $S^{n-1}_c$, as follows:
$$I=\int_{S^{n-1}_c}z_{i_1}\bar{z}_{i_2}\ldots z_{i_{2l-1}}\bar{z}_{i_{2l}}\,dz$$

Now by using the same argument as in the proof of Proposition 4.2, but this time with transformations of type $p\to\lambda p$ with $|\lambda|=1$, we see that $I$ vanishes, unless each $z_a$ appears as many times as $\bar{z}_a$ does, and this gives the first assertion.

Assume now that we are in the non-vanishing case. Then the $l_a$ copies of $z_a$ and the $l_a$ copies of $\bar{z}_a$ produce by multiplication a factor $|z_a|^{2l_a}$, so we have:
$$I=\int_{S^{n-1}_c}|z_1|^{2l_1}\ldots|z_n|^{2l_n}\,dz$$

Now by using the standard identification $S^{n-1}_c\simeq S^{2n-1}$, we get:
\begin{eqnarray*}
I
&=&\int_{S^{2n-1}}(x_1^2+y_1^2)^{l_1}\ldots(x_n^2+y_n^2)^{l_n}\,d(x,y)\\
&=&\sum_{r_1\ldots r_n}\begin{pmatrix}l_1\\ r_1\end{pmatrix}\ldots\begin{pmatrix}l_n\\ r_n\end{pmatrix}\int_{S^{2n-1}}x_1^{2l_1-2r_1}y_1^{2r_1}\ldots x_n^{2l_n-2r_n}y_n^{2r_n}\,d(x,y)
\end{eqnarray*}

By using the formula in Theorem 5.1, we get:
\begin{eqnarray*}
I
&=&\sum_{r_1\ldots r_n}\begin{pmatrix}l_1\\ r_1\end{pmatrix}\ldots\begin{pmatrix}l_n\\ r_n\end{pmatrix}\frac{(2n-1)!!(2r_1)!!\ldots(2r_n)!!(2l_1-2r_1)!!\ldots (2l_n-2r_n)!!}{(2n+2\Sigma l_i-1)!!}\\
&=&\sum_{r_1\ldots r_n}\begin{pmatrix}l_1\\ r_1\end{pmatrix}\ldots\begin{pmatrix}l_n\\ r_n\end{pmatrix}\frac{(2n-1)!(2r_1)!\ldots (2r_n)!(2l_1-2r_1)!\ldots (2l_n-2r_n)!}{(2n+\Sigma l_i-1)!r_1!\ldots r_n!(l_1-r_1)!\ldots (l_n-r_n)!}
\end{eqnarray*}

We can rewrite the sum on the right in the following way:
\begin{eqnarray*}
I
&=&\sum_{r_1\ldots r_n}\frac{l_1!\ldots l_n!(2n-1)!(2r_1)!\ldots (2r_n)!(2l_1-2r_1)!\ldots (2l_n-2r_n)!}{(2n+\Sigma l_i-1)!(r_1!\ldots r_n!(l_1-r_1)!\ldots (l_n-r_n)!)^2}\\
&=&\sum_{r_1}\begin{pmatrix}2r_1\\ r_1\end{pmatrix}\begin{pmatrix}2l_1-2r_1\\ l_1-r_1\end{pmatrix}\ldots\sum_{r_n}\begin{pmatrix}2r_n\\ r_n\end{pmatrix}\begin{pmatrix}2l_n-2r_n\\ l_n-r_n\end{pmatrix}\frac{(2n-1)!l_1!\ldots l_n!}{(2n+\Sigma l_i-1)!}
\end{eqnarray*}

The sums on the right being $4^{l_1},\ldots,4^{l_n}$, we get the formula in the statement.
\end{proof}

In the case of the free sphere, we already know from Theorem 4.4 that the standard coordinates $x_1,\ldots,x_n$ are asymptotically semicircular and free. However, the computation of their joint law for a fixed value of $n$ is a well-known open problem, of remarkable difficulty. The point is that the Gram matrix $G_{kn}$, which is nothing but Di Francesco's ``meander matrix'' in \cite{dif}, cannot be diagonalized explicitely.

The best result in this direction that is known so far is as follows.

\begin{theorem}
The moments of the free hyperspherical law are given by
$$\int_{S^{n-1}_+}x_1^{2l}\,dx=\frac{1}{(n+1)^l}\cdot\frac{q+1}{q-1}\cdot\frac{1}{l+1}\sum_{r=-l-1}^{l+1}(-1)^r\begin{pmatrix}2l+2\cr l+r+1\end{pmatrix}\frac{r}{1+q^r}$$
where $q\in [-1,0)$ is given by $q+q^{-1}=-n$.
\end{theorem}

\begin{proof}
This is proved in \cite{bcz}, the idea being that $x_1\in A_n^+$ can be modelled by a certain variable over $SU^q_2$, which can be studied by using advanced calculus methods.
\end{proof}

Our question is whether Theorem 5.1 and Theorem 5.2 have a free analogue.

\begin{question}
Does the liberated spherical integral
$$\int_{S^{n-1}_+}x_{i_1}\ldots x_{i_k}\,dx$$
appear as a ``free analogue'' of the quantities computed in Theorems 5.1 and 5.2?
\end{question}

The answer here is very unclear, even in the case where the indices $i_1,\ldots,i_k$ are all equal. In fact, the above question is probably closely related to Question 3.4.

Let us also mention that the meander determinant computed by Di Francesco in \cite{dif}, which appears as denominator of the abstract Weingarten-theoretical fraction expressing the integral in Question 5.4, is a product of Chebycheff polynomials. Our question is whether some ``magic'' simplification appears when computing the fraction.

In fact, our Questions 3.4 and 5.4 should be regarded as a slight, very speculative advance on the conceptual understanding of the various formulae in \cite{bcz}, \cite{dif}.

\section{Spectral triples}

In the reminder of this paper, our goal will be to study the ``differential structure'' of the noncommutative spheres $S^{n-1}_\times$. Besides of being of independent theoretical interest, this study will lead via the results in \cite{bve}, \cite{bg2} to a ``global look'' to our 3 spheres.

The natural framework for the study of noncommutative objects like $S^{n-1}_\times$ is Connes' noncommutative geometry \cite{con}, where the basic definition is as follows.

\begin{definition}
A compact spectral triple $(A,H,D)$ consists of the following:
\begin{enumerate}
\item $A$ is a unital $C^*$-algebra.

\item $H$ is a Hilbert space, on which $A$ acts.

\item $D$ is a (typically unbounded) self-adjoint operator on $H$, with compact resolvents, such that $[D,a]$ has a bounded extension, for any $a$ in a dense $\ast$-subalgebra (say ${\mathcal A}$) of $A$.
\end{enumerate}
\end{definition}

This definition is of course over-simplified, as to best fit with the purposes of the present paper. We refer to \cite{con} for the exact formulation of the axioms.

In what follows we will be mainly interested in the sphere $S^{n-1}$, and in its noncommutative versions $S^{n-1}_*$ and $S^{n-1}_+$. These objects are all quite simple, geometrically speaking, and we will make only a moderate use of the general machinery in \cite{con}. 

Our guiding examples, all very basic, will be as follows.

\begin{proposition}
Associated to a compact Riemannian manifold $M$ are the following spectral triples $(A,H,D)$, with $A$ and ${\mathcal A}$ being the algebra of continuous functions and that of smooth functions on $M$ respectively:
\begin{enumerate}
\item $H$ is the space of square-integrable spinors, and $D$ is the Dirac operator.

\item $H$ is the space of forms on $M$, and $D$ is the Hodge-Dirac operator $d+d^*$.

\item $H=L^2(M,dv)$, $dv$ being the Riemannian volume, and $D=\sqrt{d^*d}$. 
\end{enumerate}
\end{proposition}

Here in the first example $M$ is of course assumed to be a spin manifold. The fact that all the above triples satisfy Connes' axioms in Definition 6.1 comes from certain standard results in global differential geometry, and we refer here to \cite{con} and references therein. Let us also remark that the third example, though rather uninteresting from the viewpoint of algebraic topology or K-theory, contains all the useful information about the Riemannian geometry of the manifold, like the volume or the curvature. 

Let us go back now to our 3 noncommutative spheres, described by the algebras $A_n^\times$ in the previous sections. It is technically convenient at this point to slightly enlarge our formalism, by starting with the following ``minimal'' set of axioms.

\begin{definition}
A spherical algebra is a $C^*$-algebra $A$, given with a family of generators $x_1,\ldots,x_n$ and with a faithful positive unital trace $tr:A\to\mathbb C$, such that:
\begin{enumerate}
\item $x_1,\ldots,x_n$ are self-adjoint.

\item $x_1^2+\ldots+x_n^2=1$.

\item $tr(x_i)=0$, for any $i$.
\end{enumerate}
\end{definition}

As a first observation, each $A_n^\times$ is indeed a spherical algebra in the above sense.

We know that for $A_n=C(S^{n-1})$, there are at least 3 spectral triples that can be constructed, namely those in Proposition 6.2. In the case of $A_n^*,A_n^+$, however, or more generally in the case of an arbitrary spherical algebra, the situation with the first two constructions is quite unclear, and the third construction will be our model. 

We agree to view the identity $1$ as a length 0 word in the generators $x_1,\ldots,x_n$.

\begin{theorem}
Associated to any spherical algebra $A=<x_1,\ldots,x_n>$ is the compact spectral triple $(A,H,D)$, where the dense subalgebra ${\mathcal A}$ is the linear span of all the finite words in the generators $x_i$, and $D$ acting on $H=L^2(A,tr)$ is defined as follows:
\begin{enumerate}
\item Let $H_k=span(x_{i_1}\ldots x_{i_r}|i_1,\ldots,i_r\in\{1,\ldots,n\}, r\leq k)$.

\item Let $E_k=H_k\cap H_{k-1}^\perp$, so that $H=\oplus_{k=0}^\infty E_k$.

\item We set $Dx=kx$, for any $x\in E_k$.
\end{enumerate} 
\end{theorem}

\begin{proof}
We have to show that $[D,T_i]$ is bounded, where $T_i$ is the left multiplication by $x_i$. Since $x_i\in A$ is self-adjoint, so is the corresponding operator $T_i$. Now since $T_i(H_k)\subset H_{k+1}$, by self-adjointness we get $T_i(H_k^\perp)\subset H_{k-1}^\perp$. Thus we have: 
$$T_i(E_k)\subset E_{k-1}\oplus E_k\oplus E_{k+1}$$

This gives a decomposition of type $T_i=T_i^{-1}+T_i^0+T_i^1$. It is routine to check that we have $[D,T^\alpha_i]=\alpha T^\alpha_i$ for any $\alpha\in\{-1,0,1\}$, and this gives the result.
\end{proof}

As a first example, in the classical case the situation is as follows.

\begin{theorem}
For the algebra $A_n=C(S^{n-1})$, the spectral triple constructed in Theorem 6.4 essentially coincides with the one described in Proposition 6.2 (3). More precisely, the Dirac operator $D$ of Theorem 6.4 is related to $\sqrt{d^*d}$ by the bijective correspondence: $D=f(\sqrt{d^*d})$, where $f(s)=1-\frac{n}{2}+\frac{1}{2}\sqrt{4s^2+(n-2)^2}$, $s \in [0, \infty)$. In particular, the eigenspaces of $D$ and $\sqrt{d^*d}$ coincide. 
\end{theorem}

\begin{proof}
This follows from the well-known fact that $\sqrt{d^*d}$ diagonalizes as in Theorem 6.4, with the corresponding eigenvalues being $k(k+n-2)$, with $k=0,1,2,\ldots$
\end{proof}

\section{Quantum isometries}

We know from the previous section that associated to any spherical algebra $A$, and in particular to the algebras $A_n^\times$, is a certain spectral triple $(A,H,D)$. In the classical case $A=A_n$ this spectral triple is the one coming from the operator $D=\sqrt{d^*d}$.

Let us recall now the definition of the quantum isometry groups from \cite{bg2}, slightly modified as to fit with our setting. Let $S=(A,H,D)$ be a spectral triple of compact type, with $H$ assumed to be the GNS space of a certain faithful trace $tr:A\to\mathbb C$.

Consider the category of compact quantum groups acting on $S$ isometrically, that is, the compact quantum group (say $Q$) must have a unitary representation $U$ on $H$ which commutes with $D$, satisfies $U1_A=1_Q \otimes 1_A$ and  $ad_U$ maps $A''$ into itself. 

If this category has a universal object, then this universal object (which  is unique up to isomorphism) will be denoted by $QISO(S)$. See \cite{bg2} for more details.

\begin{proposition}
Let $A$ be a spherical algebra, and consider the associated spectral triple $S=(A,H,D)$. Then $QISO(S)$ exists.
\end{proposition}

\begin{proof}
The  proposition follows from Theorem 2.24 of \cite{bg2}, since the linear space spanned by $1_A$ is an eigenspace of $D$.
\end{proof}

\begin{theorem}
$QISO(S^{n-1}_\times)=O_n^\times$. 
\end{theorem}

\begin{proof}
Consider the standard coaction $\Phi:A_n^\times\to C(O_n^\times)\otimes A_n^\times$. This extends to a unitary representation on the GNS space $H_n^\times$, that we denote by $U$.

We have $\Phi(H_k)\subset C(O_n^\times)\otimes H_k$, which reads $U(H_k)\subset H_k$. By unitarity we get as well $U(H_k^\perp)\subset H_k^\perp$, so each $E_k$ is $U$-invariant, and $U,D$ must commute. That is, $\Phi$ is isometric with respect to $D$, and $O^\times_n$ must be a quantum subgroup of $QISO(S^{n-1}_\times)$. 

Assume now that $Q$ is compact quantum group with a unitary representation $V$ on $ H^\times$ commuting with $D$, such that $ad_V$ leaves $(A_n^\times)''$ invariant.  Since $D$ has an eigenspace consisting exactly of $x_1,\ldots,x_n$, both $V$ and $V^*$ must preserve this subspace, so we can find self-adjoint elements $b_{ij}\in C(Q)$ such that:
$$ad_V(x_i)=\sum_j b_{ij}\otimes x_j.$$ 

From the unitarity of $V$, it is also easy to see that $ad_V$ is trace-preserving, and by using this it follows that $(( b_{ij} ))$ as well as $(( b_{ji} ))$ are unitaries.  It follows in particular that the antipode $\kappa$ of $Q$ must send $b_{ij}$ to $b_{ji}$. Moreover,  using the defining relations satisfied by the $x_i$'s and the fact that $ad_V$ and $( \kappa \otimes {\rm id})\circ ad_V$ are $\ast$-homomorphism, we can prove that the $b_{ij}$'s will satisfy the same relations as those of the generators $u_{ij}$ of $C(O_n^\times)$. Indeed, for the free case there is nothing to prove, and we have verifed such relations for the classical (commutative) case, i.e. for $C(S^{n-1})$, in \cite{bg1}, the proof of which will  go through almost verbatim for the half-liberated case too, replacing the words $x_ix_j$ of length two by the length-3 words $x_ix_jx_k$.  This shows that $C(Q)$ is a quotient of $C(O_n^\times)$, so $Q$ is a quantum subgroup of $O_n^\times$, and we are done.  
\end{proof}

There are several questions raised by the above results, concerning the axiomatization of the noncommutative spheres. Perhaps the most important is the following one:

\begin{question}
What conditions on a spherical algebra $A$ ensure the fact that the corresponding quantum isometry group is ``easy'' in the sense of \cite{bsp}?
\end{question}

An answer here would of course provide an axiomatization of the ``easy spheres'', and our above results would translate into a 3-fold classification for the easy spheres, because of the classification results for easy quantum groups in \cite{bve}.

\section{Concluding remarks}

We have seen in this paper that the usual sphere $S^{n-1}$, the half-liberated sphere $S^{n-1}_*$, and the free sphere $S^{n-1}_+$, share a number of remarkable common properties. 

The general axiomatization and study of these 3 noncommutative spheres has raised a number of concrete questions, notably in connection with the general structure of the associated projective spaces (Question 3.4), with the computation of the associated spherical integrals (Question 5.4), and with the general axiomatization problem (Question 7.3). We intend to come back to these questions in some future work. 

In addition, there are many questions about what happens in the ``untwisted'' case, and in the ``non-easy'' case. Some results here are already available from \cite{bg3}.

Finally, we have the more general problem of understanding the notion of liberation and half-liberation for more general manifolds. In the 0-dimensional case it is probably possible to use the results in \cite{bgs} in order to reach to some preliminary results. In the continuous case, however, the situation so far appears to be quite unclear.

\end{document}